\documentclass[a4paper,11pt]{amsart}
\usepackage[latin1]{inputenc}
\usepackage{latexsym}
\usepackage{amssymb}
\usepackage{amsmath}
\usepackage{amsthm}
\usepackage[all]{xy}
\usepackage[dvips]{graphicx}

\newtheorem*{teos}{Theorem}

\newtheorem{teo}{Theorem}[section]
\newtheorem{lem}[teo]{Lemma}
\newtheorem{prop}[teo]{Proposition}
\newtheorem{cor}[teo]{Corollary}
\theoremstyle{definition}

\newtheorem{exe}[teo]{Example}
\newtheorem{oss}[teo]{Remark}

\newcommand{\lrg}{\longrightarrow}

\newcommand{\C}{\mathbb{C}}

\newcommand{\shA}{\mathcal{A}}

\newcommand{\shO}{\mathcal{O}}

\newcommand{\de}{\partial}
\newcommand{\debar}{\overline{\partial}}

\newcommand{\MC}{\operatorname{MC}}
\newcommand{\Def}{\operatorname{Def}}

\newcommand{\Hom}{\operatorname{Hom}}
\newcommand{\Htp}{\operatorname{Htp}}

\newcommand{\coker}{\operatorname{coker}}

\newcommand{\contr}{{\mspace{1mu}\lrcorner\mspace{1.5mu}}}

\newcommand{\bi}{\boldsymbol{i}}

\newcommand{\Art}{\mathbf{Art}}
\newcommand{\Set}{\mathbf{Set}}

\newcommand{\cil}{\operatorname{C}_{(h,g)}}
\newcommand{\hil}{\operatorname{H}_{(h,g)}}
\begin{document}

\title[Semiregularity map]{A semiregularity map annihilating
obstructions to deforming holomorphic maps}

\author{Donatella Iacono}
\address{\newline Institut f\"ur Mathematik,\hfill\newline Johannes
Gutenberg-Universit\"at, \hfill\newline Staudingerweg 9, D 55128
Mainz, Germany.} \email{iacono@uni-mainz.de}

\begin{abstract}
We study infinitesimal  deformations of holomorphic maps  of
compact, complex, K\"ahler manifolds. In particular, we describe a
generalization of Bloch's semiregularity map that annihilates
obstructions to deform   holomorphic maps  with fixed codomain.

\end{abstract}

\maketitle

\bigskip

{\bf Keywords.} Semiregularity Map, Obstruction Theory, Functors
of Artin Rings, Differential Graded Lie Algebras.

\bigskip

\section{Introduction}

The investigation of obstruction spaces plays a fundamental role
in the study of deformation theory and moduli spaces. For
instance, the obstruction theory is used to determine the
dimension of moduli spaces or the virtual fundamental class (see,
for example, \cite{bib behr-fante}, \cite{bib bloch}, \cite{bib
buchFlenner}, \cite{bib Fantec-Manet 2}, \cite{bib Kawamata}).
From the local point of view, given an infinitesimal deformation
of a geometric object, we would like to know whether it is
possible to extend this deformation or not. The idea is to
consider the same problem of extension for the associated
deformation functor. More precisely, let $F:\Art \lrg \Set$ be a
functor of Artin rings, i.e., a covariant functor from the
category $\Art$    of local Artinian $\C$-algebras (with residue
field $\C$) to the category $\Set$ of sets, such that
$F(\C)=\{\mbox{point}\}$.

A (complete) obstruction space for $F$ is a vector space $V$, such
that,  for each small extension   $0\lrg J \lrg B \lrg A\lrg 0$ in
$\Art$ and each element $x \in F(A)$, there exists an obstruction
element $v_x \in V$, associated with $x$, that is zero if and only
if $x$ can be lifted to $F(B)$.

Since this space controls the liftings, we would like to describe it, as well as possible, and know whether
the associated obstruction element is zero or not. In general, we
just know a vector space that contains the obstructions  but we have
no explicit  description of which elements of the vector space are actually
obstructions. For example, if $W$ is another vector space
which contains $V$, then also $W$ is an obstruction space for $F$.

In \cite{bib Fantec-Manet}, B.~Fantechi and M.~Manetti proved the
existence of the  \lq \lq smallest" obstruction space for
 functors associated with deformations of geometric objects.
More precisely, they proved the existence of the universal
obstruction space for deformation functors, i.e., functors of
Artin rings satisfying Schlessinger's conditions $(H_1)$ and  a
stronger version of $(H_2)$ (see \cite[Theorem~2.1]{bib Artin} and
also \cite[Lemma~ 2.11 and Theorem~6.11]{bib Fantec-Manet}).

Since it is quite difficult to determine all  the obstruction space, the idea is to  start from studying some special and easier obstructions.
In this setting, a  very useful tool is the result known as Ran-Kawamata
$T^1$-lifting theorem: if the functor is pro-representable and if
it has no  \lq \lq curvilinear obstructions'', then the functor
has no obstructions at all. Recall that the curvilinear
obstructions are the ones arising from the  curvilinear extensions
$$
0 \lrg \C  \stackrel{\cdot x^{n}}{\lrg}\frac{\C[x]}{(x^{n+1})}
\lrg \frac{\C[x]}{(x^n)} \lrg 0.
$$

This theorem  was generalized by B.~Fantechi and M.~Manetti: if $F$ is a deformation  functor, then
$F$ has no obstructions if and only if  $F$ has no curvilinear
obstructions \cite[Corollary~6.4]{bib Fantec-Manet}.

Thus, in some cases, this result guarantees that it is enough to study the curvilinear obstructions to determine   the obstruction space. More precisely, if the curvilinear obstructions  vanish,  then all the obstructions are zero.

%This result leads to the construction of the so called semi
%regularity map to reduce the obstruction space. More precisely,
%suppose that $V$ is an obstruction space for a functor $F$, then a
%semiregularity mao is a map $\sigma$ defined on $F$ that contains
%the curvilinear obstructions in the kernel. Thus, $\sigma$
%annihilates all curvilinear obstructions and so if it is
% injective, then there are no obstructions.

A fundamental fact to note is that the curvilinear obstructions  do not
generate the obstruction space. Therefore, if these obstructions   do not vanish, we have no enough information to determine the obstruction space   \cite{bib Kawamata} and \cite[Example ~5.7~(1)]{bib Fantec-Manet}.

\smallskip

In the case of infinitesimal  deformations of complex compact
manifolds,   an obstruction space  is the second cohomology vector space $H^2(X, \Theta_X)$ of the holomorphic tangent bundle $\Theta_X$ of $X$. If $X$ is also K\"{a}hler, then A.~Beauville, H.~Clemens \cite{bib
clemens} and Z.~Ran \cite{bib Ran ostruzio} \cite{bib Ran
ostruzioII} proved that the obstructions are contained in a
subspace of $H^2(X, \Theta_X)$ defined as the kernel of a well
defined map. This is the so-called \lq \lq Kodaira's principle"
(see, for example, \cite[Theorem~10.1]{bib clemens},
\cite[Corollary~3.4]{bib mane COSTRAINT},
\cite[Corollary~12.6]{bib fiorenz-Manet}, \cite[Theorem~0]{bib Ran
ostruzio} or \cite[Corollary~3.5]{bib Ran ostruzioII}).

In the case of embedded deformations of a submanifold $X$ in a
fixed manifold $Y$,  the obstructions are naturally contained in
the first cohomology vector space $H^1(X, N_{X|Y})$ of the normal
bundle $N_{X|Y}$ of $X$ in $Y$. In this case too, if $Y$ is
K\"{a}hler, then it is possible to define a map on $H^1(X,
N_{X|Y})$, called the \lq \lq semiregularity map", whose kernel contains  the curvilinear obstructions. The idea of this map
is due to S.~Bloch \cite{bib bloch}. Thus, if we can prove that this map is injective, then we conclude that the deformations are unobstructed.

Recently, M.~Manetti studied these deformations  using the
differential graded Lie algebras (DGLAs for short) and he proved
that the semiregularity map annihilates all obstructions
\cite[Theorem~0.1 and Section~9]{bib ManettiPREPRINT}.

Therefore, even if this map is not injective, we have a control on
the obstruction space, i.e., it is contained in the kernel of the
map.

\smallskip

Inspired by this work, we follow the approach via DGLA, to study the
obstructions to infinitesimal deformations of holomorphic maps of complex compact manifolds.
In \cite{bib HorikawaI},  E.~Horikawa proved that  the
obstructions to the deformations of $f:X\lrg Y$, with fixed
codomain, are contained in the second cohomology vector space $
H^2(C^\cdot_{f_*}) $ of the cone $C^\cdot_{f_*}$, associated with
the complex morphism $f_*: A_X^{0,*}(\Theta_X) \lrg
A_X^{0,*}(f^*\Theta_Y).$

\smallskip

Using the approach via DGLAs, we can give an easy proof of this
theorem (Proposition~\ref{prop mia uguale horikawa}) and,
furthermore,   we can improve it in the case of K\"{a}hler
manifolds. Our main result is the following theorem
(Corollary~\ref{coroll semiregolari-ostruzion}).

\begin{teos}
Let $f:X \lrg Y$ be a holomorphic map of  compact K\"{a}hler
manifolds. Let $p=dim Y - dim X$. Then, the obstruction space to
the infinitesimal deformations of $f$  with fixed $Y$  is
contained in the kernel of the map
$$
\sigma :  H^2(C^\cdot_{f_*}) \lrg H^{p+1}(Y,\Omega_Y^{p-1}).
$$
\end{teos}

In the case of an inclusion $X \hookrightarrow Y$, the previous
map reduces to Bloch's semiregularity map. We remark that this map
annihilates all obstructions.

In \cite{bib buchFlenner}, R.-O.~Buchweitz and H.~Flenner studied
deformations of coherent modules and, as a particular case,
deformations of holomorphic maps. They used very different
techniques and they also produced a semiregularity map
\cite[Theorem~7.23]{bib buchFlenner}, but they didn't explicitly
state that their  map annihilates all obstructions (and not merely
the curvilinear ones).

\bigskip

{ \bf Acknowledgments.} It is a pleasure for me to show my deep
gratitude to the advisor of my PhD thesis Prof. Marco Manetti, an
excellent   helpful professor who supports and encourages me every
time. I'm indebted with him for many useful discussions, advices
and suggestions. Several ideas of this work are grown under his
influence.  The first draft of this paper was written at the
Mittag-Leffler Institute in Stockholm, during the special year on
Moduli Spaces. I am very grateful for the support received. I am
lalso grateful to  the referee for improvements in the exposition
of the paper.

\section{Notation}

We will work on the field of complex number  $\C$ and all vector
spaces, linear maps, tensor products etc. are intended over $\C$.

If $A$ is an object in $\Art$, then $m_A$ denotes its maximal
ideal.

Unless otherwise specified,  by a manifold we mean a compact,
(complex) connected and smooth variety.

Given a manifold $X$, we denote by $\Theta_X $ the holomorphic
tangent bundle, by $\shA_X^{p,q}$ the sheaf of differentiable
$(p,q)$-forms on $X$ and by $A_X^{p,q}=\Gamma(X,\shA_X^{p,q})$ the
vector space of  its global sections. More generally,
$\shA_X^{p,q}(\Theta_X)$ is the sheaf of differentiable
$(p,q)$-forms on $X$ with values in  $\Theta_X $ and
$A_X^{p,q}(\Theta_X)=\Gamma(X,\shA_X^{p,q}(\Theta_X))$ is the
vector space of  its global sections.

Finally,  by a map $f:X \lrg Y$ we always mean  a holomorphic
morphism of (complex compact) manifolds and  we denote by  $f^*$
and $f_*$ the induced maps, i.e.,
$$
f^*: A_Y^{p,q}(\Theta_Y) \lrg A_X^{p,q}(f^*\Theta_Y) \quad \mbox{
and} \quad f_*: A_X^{p,q}(\Theta_X) \lrg A_X^{p,q}(f^*\Theta_Y).
$$
The cone $C^\cdot_{f_*}$
is  the complex
$(C^\cdot_{f_*},D)$ with

$$ C^i_{f_*}:=
A_X^{0,i}(T_X) \oplus A_X^{0,i-1}(f^*T_Y)
$$
and
$$
 D:C^i_{f_*} \lrg C^{i+1}_{f_*},
$$
$$
(l,n)\mapsto (\debar l, f_*(l) -\debar n) \in A_X^{0,i+1}(T_X) \oplus
A_X^{0,i}(f^*T_Y).
$$

\section{The semiregularity map}\label{sezio def semiregular map s}

Let $f:X \lrg Y$ be a  map of   K\"{a}hler manifolds,  $n=dim X$
and  $p=dim Y - dim X$. Let    $\mathcal{H}$ be the space of
harmonic forms on $Y$ of type $(n+1,n-1)$. By Dolbeault's theorem
and Serre's duality, we have $\mathcal{H}^\nu=(H^{n-1}
(Y,\Omega_Y^{n+1}))^\nu= H^{p+1}(Y,\Omega_Y^{p-1})$.

Using the contraction $\contr$ of vector fields with differential
forms, for each $\omega \in \mathcal{H}$, we can define the
following map
$$
A_X^{0,*}(f^*\Theta_Y) \stackrel{\contr \omega}{\lrg}
A_X^{n,*+n-1},
$$
$$
\contr \omega( \phi f^* \chi)=\phi f^*(\chi \contr \omega)\in
A_X^{n,p+n-1}  \qquad \forall\  \phi f^* \chi \in
A_X^{0,p}(f^*\Theta_Y).
$$
It can be proved (see Lemma~\ref{lemma f^*(x-|w)=f_*x-|f^*w}) that
if $f^* \omega=0$, then the following diagram
\begin{center}
$\xymatrix{A_X^{0,*}(f^*\Theta_Y) \ar[r]^{\ \ \contr \omega}
& A_X^{n,*+n-1} \\
A_X^{0,*}(\Theta_X) \ar[u]^{f_*} \ar[r]  & 0 \ar[u]     \\
}$
\end{center}
is commutative. Thus, for each $\omega$, we get  a morphism
$$
H^2(C^\cdot_{f_*}) \lrg
 H^{n}(X,\Omega^n_ X),
$$
which, composed with the integration on $X$, gives the
semiregularity map
$$
\sigma :  H^2(C^\cdot_{f_*})  \lrg H^{p+1}(Y,\Omega_Y^{p-1}).
$$

If $f$ is the inclusion map $X\hookrightarrow Y$, then $
H^2(C^\cdot_{f_*}) \cong H^1(X,N_{X|Y})$,  where $N_{X|Y}$ is the
normal bundle of $X$ in $Y$. In this case, the previous map
$\sigma$ reduces to Bloch's semiregularity map (see \cite{bib
bloch} or \cite[Section~9]{bib ManettiPREPRINT}), i.e.,
$$
\sigma: H^1(X,N_{X|Y}) \lrg H^{p+1}(Y,\Omega_Y^{p-1}).
$$

\begin{exe}

Let $S$ be a K3 surface. Then, the canonical bundle is trivial,
$\Theta_S \cong \Omega^1_S$, $q(s)=dim H^1(S,\shO_S)=0$ and
$p_g(S)= dim H^2(S,\shO_S)=1$.

Let $f:C \lrg S$ be a non constant holomorphic map  from a smooth
curve $C$ in $S$ (the differential $f_*:\Theta_C \lrg f^*\Theta_S$
is non zero at the generic point). If we consider  the deformation
of $f$, with fixed codomain $S$, then  the semiregularity map
$$
\sigma:H^2(C^\cdot_{f_*})  \lrg H^{2}(S,\shO_S)\cong \C
$$
is surjective. Indeed, let  $N_f$ be the cokernel of $f_*$, i.e.,
$$
\Theta_C \stackrel{f_*}{\lrg} f^*\Theta_S\lrg N_f\lrg 0.
$$
The hypothesis on $f_*$ implies that the sequence
$$
0\lrg \Theta_C \stackrel{f_*}{\lrg} f^*\Theta_S\lrg N_f\lrg 0
$$
is also exact. Therefore, $H^i(C,N_f)\cong H^{i+1}
(C^\cdot_{f_*})$, for each $i \geq 0$, and   the induced map
$$
H^1(C, f^*\Theta_S)  \lrg H^1(C,N_f)
$$
is surjective.

Consider the pull-back $f^* \Omega^1_S \lrg \Omega ^1_C$ and
denote by $\mathcal{L}$ and $\Delta$ the kernel and the cokernel,
respectively, i.e.,
\begin{center}
$\xymatrix{0\ar[r] & \mathcal{L}   \ar[r] &  f^* \Omega^1_S
\ar[rd] \ar[rr]
& &\Omega ^1_C  \ar[r]&  \Delta \ar[r] & 0 \\
 & & & K  \ar[ru]\ar[rd]& & &       \\
 & & 0 \ar[ru]& & 0.& &       \\
}$
\end{center}
By hypothesis on $f$, $\Delta$ is a torsion sheaf and so
$H^1(C,\Delta)=0$. Therefore,  $H^1(C,K) \lrg H^1(C,\Omega^1_C)$
is surjective.

Moreover, $H^2(C,\mathcal{L})=0$ and so $H^1(C,f^*\Omega^1_S)\lrg
H^1(K)$ is surjective. In conclusion, the induced map
$$
H^1(C,f^*\Omega^1_S) \lrg  H^1(C,\Omega^1_C)
$$
is surjective. By the integration on $C$, we get a surjective map
$$
H^1(C,f^*\Omega^1_S) \lrg  \C.
$$
Since the diagram
\begin{center}
$\xymatrix{ &  H^1(C,f^* \Omega^1_S) \ar@{>>}[ld] \ar@{>>}[rd] &\\
  H^1(C,N_f)\ar[rr]^{\sigma}     &   & \C   \\
}$
\end{center}
is commutative, the semiregularity map is surjective.

\end{exe}

%In particular, if $C\subset S$ is an hyperplane section, S.~Bloch
%proved \cite[Proposition~1.1]{bib bloch} that the semiregularity
%map
%$$
%H^1(N_{C|S})\lrg H^2(S, \shO_S )
%$$
%arises as the boundary map in the cohomology sequence associated
%to
%$$
%0\lrg \shO_S \lrg \shO_S (C) \lrg N_{C|S} \lrg 0.
%$$
%Since  $H^1(S, \shO_S (C))=H^2(S, \shO_S (C))=0$, the
%semiregularity map is an isomorphism. Thus, we recover the well
%known fact that the   hyperplane sections are unobstructed.

\section{Proof of the main Theorem}

Nowadays, the approach to deformation theory via DGLAs is quite
standard (see for example \cite{bib buchMilso}, \cite{bib
kontsevich}, \cite{bib manRENDICONTi}).

In \cite{bib ManettiPREPRINT}, M.~Manetti used the DGLAs to study
the obstructions of the inclusion map and   Bloch's semiregularity
map.

Inspired by his work, we also prove our main theorem using the
DGLAs and, in  particular,   the techniques developed in \cite{bib
tesidottorato} and \cite{bib articolo1tesi}.

For reader's convenience, we recall the main results of these
papers.

To study   deformations of holomorphic maps via DGLAs, it is
convenient to define a deformation functor associated with a pair
of morphisms of DGLAs. More precisely, let  $h :L \lrg M$ and $g:
N \lrg M$ be  morphisms of DGLAs, with $M$ concentrated in non
negative degrees, i.e.,
\begin{center}
$\xymatrix{ & L \ar[d]^h   \\
N\ar[r]^{g}   & M. \\ }$
\end{center}
Then, the deformation functor associated with the pair $(h,g)$ is
$$
\Def_{(h,g)}: \Art \lrg \Set,
$$
$$
\Def_{(h,g)}(A) =\frac{\MC_{(h,g)}(A)}{gauge},
$$
where
$$
\MC_{(h,g)}(A)= \{(x,y,e^p) \in (L^1 \otimes m_ A)\times  (N^1
\otimes m_A ) \times \operatorname{exp}(M^0 \otimes m_A)  |
$$
$$
dx + \displaystyle\frac{1}{2} [x,x]=0,\   dy +\frac{1}{2}[y,y]=0,\
g(y)=e^p*h(x) \},
$$
and the gauge equivalence is induced by the   gauge action of
$\exp(L^0 \otimes m_A) \times \exp(N^0 \otimes m_A)$ on
$\MC_{(h,g)}(A)$, given by
$$
(e^a,e^b)*(x,y,e^p)= (e^a*x,e^b*y,e^{g(b)}e^pe^{-h(a)}).
$$
Let $(\cil^\cdot,D)$  be the differential graded vector space with
$$
\cil^i=L^i \oplus N^i \oplus M^{i-1} \quad \mbox{ and } \quad
D(l,n,m)=( dl, d n,-d m-g(n)+h(l)).
$$
Then, the tangent space of $\Def_{(h,g)}$ is $H^1 (\cil^\cdot)$
and the obstruction space of $\Def_{(h,g)}$ is naturally contained
in $H^2 (\cil^\cdot)$ \cite[Lemma~III.1.19]{bib tesidottorato} or
\cite[Section~4.2]{bib articolo1tesi}.

\begin{lem}\label{lemma L M N abelian ->DEF(h,g)liscio}
Let $h:L \lrg M$ and $g:N \lrg M$ be morphisms of abelian DGLAs.
Then, the functor $\MC_{(h,g)}$ is smooth, that is, it has no
obstructions.
\end{lem}
\begin{proof}
See \cite[Lemma~II.1.20]{bib tesidottorato}.

\end{proof}

\begin{oss}
Every commutative diagram of morphisms of DGLAs
\begin{center}
$\xymatrix{ & L\ar[d]_h
\ar[r]^{\alpha'} & P\ar[d]^\eta \\
& M\ar[r]^\alpha & Q \\
 N\ar[ur]^g \ar[r]^{\alpha''} & R \ar[ur]_\mu & \\}$
\end{center}
induces a \emph{morphism $\varphi^\cdot$ of complexes }
$$
\cil^i \ni (l,n,m) \stackrel{\varphi ^i}{\longmapsto} (\alpha'(l),
\alpha''(n),\alpha(m)) \in \operatorname{C}^i_{(\eta,\mu)}
$$
and a natural \emph{transformation $F$ }of the associated
deformation functors, i.e.,
$$
F: \Def_{(h,g)} \lrg \Def_{(\eta,\mu)}.
$$

\end{oss}

\begin{prop}
\label{prop NO exte quasi iso C(h,g)-C(n,m)then DEf iso DEF} If
$\varphi^\cdot:\cil^\cdot \lrg
\operatorname{C}^\cdot_{(\eta,\mu)}$ is a quasi-isomorphism of
complexes, then $F:\Def_{(h,g)} \lrg \Def_{(\eta,\mu)}$ is an
isomorphism of functors.
\end{prop}
\begin{proof}
See \cite[Theorem~III.1.23]{bib tesidottorato}.
\end{proof}

\begin{prop}\label{prop DEF_>DEF liscio,ostruz in ker tra H^2}
Let
\begin{center}
$\xymatrix{ & L\ar[d]_h
\ar[r]^{\alpha'} & P\ar[d]^\eta \\
& M\ar[r]^\alpha & Q \\
N\ar[ur]^g \ar[r]^{\alpha''} & R \ar[ur]_\mu & \\}$
\end{center}
be a commutative diagram of differential graded Lie algebras. If
the functor $\Def_{(\eta,\mu)}$ is smooth, then the obstruction
space of $\, \Def_{(h,g)}$ is contained in the kernel of the map
$$
H^2(\cil^\cdot) \lrg H^2(\operatorname{C}_{(\eta,\mu)}^\cdot).
$$
\end{prop}
\begin{proof}
The natural transformation $F:\Def_{(h,g)} \lrg \Def_{(\eta,\mu)}$
induces a linear map between the obstruction spaces. If
$\Def_{(\eta,\mu)}$ is smooth, then its obstruction space is zero.

\end{proof}

By a suitable choice of the morphisms $h$ and $g$, we can study
the infinitesimal deformations of holomorphic maps.

Indeed, let $f:X \lrg Y$ be a holomorphic map, $Z=X\times Y$ and
$\Gamma\subset Z$ the graph of $f$. Let
$$
F:X \lrg \Gamma \subseteq Z:=X \times Y,
$$
$$
\ \ \ x \longmapsto (x,f(x)),
$$
and $p:Z\lrg X$ and $q:Z \lrg Y$ be the natural projections.

Then, we have the following commutative diagram:
\begin{center}
$\xymatrix{X \ar[rrrr]^{F} \ar[ddrrr]_{id} \ar[ddrrrrr]^f & &  & &
Z\ar[ddl]_{ p} \ar[ddr]_{ \ \ q} & \\
& & &  & &  \\
& & & X & & Y. \\ }$
\end{center}
In particular, $F^*\circ p^*=id$ and $F^* \circ q^*=f^*$. Since
$\Theta_Z=p^*\Theta_X \oplus q^*\Theta_Y$, it follows that
$F^*(\Theta_{Z})=\Theta_X \oplus f^*\Theta_Y$. Define the morphism
$\gamma:\Theta_Z \lrg f^*\Theta_Y$ as the product
$$
\gamma: \Theta_Z \stackrel{F^*}{\lrg}\Theta_X \oplus f^* \Theta_Y
\stackrel{(f_*,-id)}{\lrg} f^*\Theta_Y;
$$
moreover, let $\pi$ be the following surjective morphism:
$$
A^{0,*}_{Z}(\Theta_{Z}) \stackrel{\pi}{\lrg}
A_X^{0,*}(f^*\Theta_Y) \lrg 0,
$$
$$
\pi(\omega\, u)= F^*(\omega)\gamma(u), \qquad \forall \ \omega \in
A^{0,*}_{Z}, \ u \in \Theta_{Z}.
$$
Since each $u \in \Theta_{Z}$ can be written as $u=p^*v_1+q^*v_2$,
for some $v_1 \in\Theta_X$ and $v_2 \in \Theta_Y$, we also have
$$
\pi (\omega u)=F^*(\omega)(f_*(v_1)-f^*(v_2)).
$$
The algebra $A^{0,*}_{Z}(\Theta_{Z})$ is the Kodaira-Spencer
(differential graded Lie) algebra of    $Z$ and we denote by
$A_{Z}^{0,*}(\Theta_{Z}(-log\, \Gamma))$ its differential graded
Lie subalgebra  defined by the following exact sequence
\begin{equation}\label{equa art2 defi L=T(log)}
0 \lrg A_{Z}^{0,*}(\Theta_{Z}(-log\, \Gamma)) \lrg
A_{Z}^{0,*}(\Theta_{Z})
  \stackrel{\pi}{\lrg} A_X^{0,*}(f^*\Theta_Y)\lrg 0.
\end{equation}

The DGLA $A_{Z}^{0,*}(\Theta_{Z})$ controls the infinitesimal
deformations of $Z$ and $A_{Z}^{0,*}(\Theta_{Z}(-log\, \Gamma))$
controls the infinitesimal deformations of the pair $\Gamma\subset
Z$, i.e., each solution of the Maurer-Cartan equation in
$A_{Z}^{0,*}(\Theta_{Z}(-log\, \Gamma))$ defines a deformation of
both $\Gamma$ and $Z$ \cite{bib ManettiPREPRINT}.

Consider the morphism of DGLAs $g=(p^*,q^*):A^{0,*}_{X }(\Theta_X)
\times A^{0,*}_{Y}(\Theta_Y) \lrg A_{Z}^{0,*}(\Theta_{Z})$. The
solutions  $n=(n_1,n_2)$ of the Maurer-Cartan equation in
$N=A^{0,*}_{X  }(\Theta_X) \times A^{0,*}_{Y}(\Theta_Y) $
correspond to  infinitesimal deformations of both $X$ (induced by
$n_1$) and  $Y$ (induced by $n_2$). Moreover, the image $g(n)$
satisfies the Maurer-Cartan equation in
$M=A^{0,*}_{Z}(\Theta_{Z})$ and so it is associated with an
infinitesimal deformation of $ Z$, that is exactly the one
obtained as product of the deformations of $X$ (induced by $n_1$)
and of $Y$ (induced by $n_2$).

Next, fix $M=A_{Z}^{0,*}(\Theta_{Z}) $,
$L=A_{Z}^{0,*}(\Theta_{Z}(-log\, \Gamma)) $, $h$ the inclusion
$L\hookrightarrow  M$, $N=A^{0,*}_{X }(\Theta_X) \times
A^{0,*}_{Y}(\Theta_Y) $ and $g=(p^*,q^*):N \lrg M$, i.e.,
\begin{equation}\label{diagram L=A M=A N=AxA in background
semiregolar} \xymatrix{ & & A_{Z}^{0,*}(\Theta_{Z}(-log\, \Gamma))
\ar@{^{(}->}[d]^{\ h}   & \\
A_X^{0,*}(\Theta_X) \times A_Y^{0,*}(\Theta_Y)
\ar[rr]^{g=(p^*,q^*)}
 & & A_{Z}^{0,*}(\Theta_{Z}).  &
 \\ }
\end{equation}
If $\Def(f)$ is  the   functor  of the infinitesimal deformations
of the map $f$, then   the following theorem holds.

\begin{teo}\label{teo def(h,g)=def(f) articolo2 semir}
Let $f:X \lrg Y$ be a holomorphic map of compact complex manifold.
Then, with the  notation above, there exists an isomorphism of
functors
$$
\Def_{(h,g)}\cong \Def (f).
$$
\end{teo}
\begin{proof}
See \cite[Theorem~IV.2.5]{bib tesidottorato} or
\cite[Theorem~5.11]{bib articolo1tesi}.
\end{proof}

Furthermore,  for each choice of the pair $(h,g)$, there exist a
DGLA  $H_{(h,g)}$ and an isomorphism $\Def_{\hil} \cong
\Def_{(h,g)}$ \cite[Corollary~6.18]{bib articolo1tesi}. In
particular, there exists an explicit description of a DGLA
$H_{(h,g)}$  that controls the infinitesimal deformations of $f$,
i.e., $ \Def(f) \cong \Def_{H(h,g)}$ \cite[Theorem~6.19]{bib
articolo1tesi}.

In general, it is not easy to handle the DGLA  $H_{(h,g)}$ and so
it is convenient to use the  functor $\Def_{(h,g)}$,  associated
with the previous  diagram (\ref{diagram L=A M=A N=AxA in
background semiregolar}).

Indeed, for example, if we want to study the infinitesimal
deformations of $f$ with fixed domain, it   suffices to take
$N=A^{0,*}_{Y }(\Theta_Y)$.

Analogously, in the case of deformations of a map $f$ with fixed
codomain $Y$, the DGLA $N$ reduces to $A_X^{0,*}(\Theta_X)$ and so
diagram (\ref{diagram L=A M=A N=AxA in background semiregolar})
reduces to

\begin{center}
$\xymatrix{ & L \ar[d]^h   & \\
A_X^{0,*}(\Theta_X) \ar[r]^{p^*}  \ar@/^/[rrd]_{f_*} &
A_{Z}^{0,*}(\Theta_{Z})\ar[dr]^\pi &    \\
& & A_X^{0,*}(f^*\Theta_Y), &
 \\ }$
\end{center}
where $f_*$ is the product $\pi \circ p^*$.

Using this diagram and Theorem~\ref{teo def(h,g)=def(f) articolo2
semir}, we can easily prove the following proposition due to
E.~Horikawa \cite{bib HorikawaI}.

\begin{prop}\label{prop mia uguale horikawa}
The tangent space to the infinitesimal deformations   of a
holomorphic map $f:X \lrg Y$, with  fixed codomain $Y$, is
$H^1(C^\cdot_{f_*})$ and the obstruction space is naturally
contained in $H^2(C^\cdot_{f_*})$.
\end{prop}
\begin{proof}
Theorem~\ref{teo def(h,g)=def(f) articolo2 semir} implies that the
infinitesimal deformation  functor of $f$,  with $Y$ fixed, is
isomorphic to $\Def_{(h,p^*)}$. Therefore, the tangent space is
$H^1(C^\cdot_{(h,p^*)})$ and the obstruction space is naturally
contained in $H^2(C^\cdot_{(h,p^*)})$. Since $h$ is injective, we
have isomorphisms  $H^i (C^\cdot_{(h,p^*)})\cong H^i (C^\cdot_{\pi
\circ p^*})=H^i(C^\cdot_{f_*})$, for each $i$.
\end{proof}

\medskip

Our main theorem improves this result in the case of K\"{a}hler
manifolds. To prove it, we need some preliminary lemmas.

\begin{lem}\label{lemma f^*(x-|w)=f_*x-|f^*w}
Let $f:X \lrg Y$ be a  holomorphic map of complex compact
manifolds. Let $\chi \in \shA_Y^{0,*}(\Theta_Y)$  and $\eta \in
\shA_X^{0,*}(\Theta_X)$ such that $f^*\chi =f_* \eta \in
\shA_X^{0,*}(f^*\Theta_Y)$. Then, for any  $\omega \in
\shA_Y^{*,*}$
$$
f^*(\chi \contr \omega)= \eta \contr f^* \omega.
$$
\end{lem}

\begin{proof}
See \cite[Lemma~II.6.1]{bib tesidottorato}. It follows from an
easy calculation in local holomorphic coordinates.
\end{proof}

\noindent Let $f:X \lrg Y$ be a holomorphic map, $Z=X\times Y$ and
$\Gamma\subset Z$ the graph of $f$.

\begin{lem}\label{lemma kaler aciclico  q^* A_Y}
If X and Y are compact $K\ddot{a}hler$ manifolds, then the
sub-complexes $Im(\de)= \de A_Z $, $\de A_\Gamma $, $\de A_{Z}
\cap q^* A_Y$ and $\de A_{Z} \cap p^* A_X$ are acyclic.
% ??????? (with respect to $\debar$)

\end{lem}
\begin{proof}
See \cite[Lemma~II.2.2]{bib tesidottorato}. It follows from the
$\de \debar$-Lemma.
\end{proof}

\begin{oss}\label{oss de-debar lemma per X,Y,XxY, Gamma}
In the previous lemma, the  K\"{a}hler hypothesis on $X$ and $Y$
can be substituted by the validity of the $\de \debar$-lemma  in
$A_X$, $A_Y$, $A_Z=A_{X \times Y}$ and $A_\Gamma$. In particular,
 it holds for every compact complex manifolds  bimeromorphic to a
K\"{a}hler manifolds  \cite[Corollary 5.23]{bib DeligneGMS}.

\end{oss}

Let $W$ be a  manifold and   $A_W^{0,*}(\Theta_W)$  its
Kodaira-Spencer algebra. Then, we define the contraction map $\bi$
as follows:
$$
\bi : A_W^{0,*}(\Theta_W) \lrg \operatorname \Hom^*(A_W,A_W),
$$
$$
\bi_a(\omega)=a \contr \omega, \qquad \forall\  a \in
A_W^{0,*}(\Theta_W) \mbox{ and }  \omega \in  A_W^{*,*}.
$$

\noindent Therefore, $\bi(A_W^{0,j}(\Theta_W))\subset \oplus_{h,l}
\Hom^0(A_W^{h,l},A_W^{h-1,l+j}) \subset \Hom^{j-1}(A_W,A_W)$.

\medskip

In order to interpret $\bi$ as a morphism of DGLAs, the key idea,
due to M.~Manetti \cite[Section~8]{bib ManettiPREPRINT}, is to
substitute $\Hom^*(A_W,A_W)$ with the differential  graded  vector
space $\Htp \left (\ker(\de), \dfrac{A_W} {\de
A_W}\right)=\bigoplus_i \Hom^{i-1} \left (\ker(\de), \dfrac{A_W}
{\de A_W}\right)$. Consider on $\Htp \left (\ker(\de), \dfrac{A_W}
{\de A_W}\right)$  the following differential $\delta$ and bracket
$\{\ ,\ \}$:
$$
\delta(f)=-\debar f - (-1)^{\deg(f)}f \debar,
$$
$$
\{f,g\}=f \de g -(-1)^{\deg(f)\deg(g)}g \de f.
$$

\begin{lem}
$\Htp \left (\ker(\de), \dfrac{A_W} {\de A_W}\right)$ is a DGLA
and the linear map
$$
\bi: A_W^{0,*}(T_W) \lrg \Htp\left(\ker(\de),\dfrac{A_W}{\de
A_W}\right)
$$
is a morphism of DGLAs.
\end{lem}
\begin{proof}
See \cite[Proposition~8.1]{bib ManettiPREPRINT}.

\end{proof}

\begin{oss}
For any pair of graded vector spaces $V$ and $W$, there exists an
isomorphism $H^i(\Htp(V,W))\cong \Htp^i(H^*(V),H^*(W)))$, for each
$i$.
\end{oss}

Next, we   apply this construction to $Z=X\times Y$. \\
Let $\Gamma$ be the graph of $f$ in $Z$ and $I_\Gamma \subset A_Z$
the space of the differential forms vanishing on $\Gamma$. The
DGLA $L= A_{Z}^{0,*}(\Theta_{Z}(-log\, \Gamma))$ defined in
(\ref{equa art2 defi L=T(log)}) satisfies the following property
$$
L \subset \{a \in  A^{0,*}_Z(\Theta_Z)| \ \bi_a(I_\Gamma)\subset
I_\Gamma\}.
$$
Furthermore,
$$
p^*A_X^{0,*}(\Theta_X)\subset \{a \in  A^{0,*}_Z(\Theta_Z)| \
\bi_a(q^*A_Y)=0\},
$$
where $p$ and $q$ are the projections of $Z$ onto $X$ and $Y$,
respectively.

In conclusion, we can define the following commutative diagram of
morphisms of DGLAs

\begin{center}
\begin{equation}\label{equa diagramma tre piani}
 \xymatrix{ L\ar@{^{(}->}[d]^h \ar[r] & K= \left\{f\in
\Htp\left(\ker(\de),\dfrac{A_Z}{\de A_Z}\right)\mid
f(I_\Gamma\cap\ker(\de))\subset \dfrac{I_\Gamma}{I_\Gamma\cap\de
A_Z}\right\}
\ar@{^{(}->}[d]^\eta \\
 A^{0,*}_Z(\Theta_Z)\ar[r]& \Htp\left(\ker(\de),\dfrac{A_Z}{\de
A_Z}\right) \\
A^{0,*}_X(\Theta_X)\ar[u]_{p^*} \ar[r] & J=\left\{f\in
\Htp\left(\ker(\de),\dfrac{A_Z}{\de A_Z}\right)\mid
f(\ker(\de)\cap q^* A_Y)=0 \right\},  \ar@{^{(}->}[u]_\mu \\ }
\end{equation}
\end{center}
where the horizontal maps are all given by $\bi$.

We note that diagram (\ref{equa diagramma tre piani}) induces a
natural transformation of deformation functors:
$$
\mathcal{I}:\Def_{(h,p^*)}\lrg \Def_{(\eta,\mu)}.
$$

\begin{lem}\label{lemma ostruzione nel ker tra H^2}
If the differential graded vector spaces $(\de A_Z,\debar)$, $(\de
A_\Gamma,\debar)$ and $(\de A_Z \cap q^*A_Y, \debar)$ are acyclic,
then the functor $\Def_{(\eta,\mu)}$ has no obstructions. In
particular, the obstruction space of $\Def_{(h,p^*)}$ is naturally
contained in the kernel of the map
$$
H^2(C^\cdot_{(h,p^*)}) \stackrel{\mathcal{I}}{\lrg}
H^2(C^\cdot_{(\eta,\mu)}).
$$
\end{lem}
\begin{proof}This lemma is an extension   of
\cite[Lemma~8.2]{bib ManettiPREPRINT}.

The projection $\ker(\de)\to \ker(\de)/\de A_Z$ induces a
commutative diagram
\begin{center}
\begin{equation}\label{eqa diagr 2 colon DEF_(n,m) abeliano}
\xymatrix{ K\ar[d]^\eta & \{f\in K | f(\de A_Z)=0\}\ar[d]_{\eta'}
\ar[l]^\alpha  \\
\Htp\left(\ker(\de),\dfrac{A_Z}{\de A_Z}\right)  &
\Htp\left(\dfrac{\ker(\de)}{\de
A_Z},\dfrac{A_Z}{\de A_Z}\right)\ar[l]^\beta\\
J \ar[u]_\mu & \{f\in J| f(\de A_Z)=0\}.\ar[u]^{\mu'}
\ar[l]^\gamma
\\ }
\end{equation}
\end{center}

Since $\de A_Z$ is acyclic, $\beta$ is a quasi-isomorphism of
DGLAs. Since
$$
\coker(\alpha)=\{ f \in \Htp\left(\de A_Z,\dfrac{A_Z}{\de
A_Z}\right) | f(I_\Gamma \cap \de A_Z)\subset
\dfrac{I_\Gamma}{I_\Gamma\cap\de A_Z}  \},
$$
there exists an exact sequence
$$
0\to \Htp\left(\dfrac{\de A_Z}{I_\Gamma\cap \de
A_Z},\dfrac{A_Z}{\de A_Z}\right) \to \coker(\alpha)\to
\Htp\left(I_\Gamma\cap \de A_Z,\dfrac{I_\Gamma}{I_\Gamma\cap \de
A_Z}\right)\to 0.
$$
Moreover, the exact sequence
$$
0 \lrg I_\Gamma \cap A_Z \lrg \de A_Z \lrg \de A_\Gamma\lrg 0
$$
implies that $I_\Gamma \cap A_Z$ and $\dfrac{\de A_Z}{I_\Gamma\cap
\de A_Z}=\de A_\Gamma$ are acyclic. Thus, the complexes
$\Htp\left(\dfrac{\de A_Z}{I_\Gamma\cap \de A_Z},\dfrac{A_Z}{\de
A_Z}\right)$ and $\Htp\left(I_\Gamma\cap \de
A_Z,\dfrac{I_\Gamma}{I_\Gamma\cap \de A_Z}\right)$ are acyclic and
the same holds for $\coker(\alpha)$, i.e., $\alpha$ is a
quasi-isomorphism.

As to  $\gamma$, we have
$$
\coker(\gamma)=
$$
$$
\{ f \in \Htp\left(\de A_Z,\dfrac{A_Z}{\de A_Z}\right) |\  f(\de
A_Z \cap q^*A_Y)=0 \}=
$$
$$
 \Htp\left(\dfrac{\de A_Z}{\de A_Z \cap
q^*A_Y},\dfrac{A_Z}{\de A_Z}\right).
$$
By hypothesis, $\de A_Z \cap q^*A_Y $ and $\de A_Z$ are acyclic
and so the same holds for $\dfrac{\de A_Z}{\de A_Z \cap q^*A_Y}$.
Then, $\coker(\gamma)$ is acyclic, i.e.,   $\gamma$ is also a
quasi-isomorphism.

Therefore, by Lemma~\ref{prop NO exte quasi iso C(h,g)-C(n,m)then
DEf iso DEF},  there exists   an isomorphism of deformation
functors $\Def_{(\eta,\mu)}\cong \Def_{(\eta',\mu')}$. We note
that the elements of the three algebras of the right column of
(\ref{eqa diagr 2 colon DEF_(n,m) abeliano}) vanish on $\de A_Z$.
Then, by the definition  of the bracket $\{ \ , \ \}$, these
algebras are abelian and so, by Lemma~\ref{lemma L M N abelian
->DEF(h,g)liscio},  the functor $\Def_{(\eta,\mu)} \cong
\Def_{(\eta',\mu')}$ has no obstructions.

Therefore, by Proposition~\ref{prop DEF_>DEF liscio,ostruz in ker
tra H^2}  the obstruction space of $\Def_{(h,p^*)}$ lies in the
kernel of $H^2(C^\cdot_{(h,p^*)}) \stackrel{\mathcal{I}}{\lrg}
H^2(C^\cdot_{(\eta,\mu)}) $.
\end{proof}

\begin{teo}\label{teo ostruzi in ker H^2}
Let $f:X \lrg Y$ be a  holomorphic map of compact K\"{a}hler
manifolds. Then, the obstruction space to the infinitesimal
deformations of $f$ with fixed codomain is contained in the kernel
of the following map
$$
H^2(C^\cdot_{f_*})\stackrel{\mathcal{J}}{\lrg}
H^1\left(\Htp(I_\Gamma \cap \ker(\de)\cap q^*A_Y,A_\Gamma)\right).
$$
\end{teo}
\begin{proof}
By Lemma~\ref{lemma kaler aciclico  q^* A_Y},   the complexes
$(\de A_Z,\debar)$, $(\de A_\Gamma,\debar)$ and $(\de A_Z \cap
q^*A_Y, \debar)$ are acyclic. Then, Lemma~\ref{lemma ostruzione
nel ker tra H^2} implies that  the obstruction space lies in the
kernel of the following map
$$
H^2(C^\cdot_{(h,p^*)}) \stackrel{\mathcal{I }}{\lrg}
H^2(C^\cdot_{(\eta,\mu)}).
$$
Since $h$ is injective, as in Proposition~\ref{prop mia uguale
horikawa}, we have $H^2(C^\cdot_{(h,p^*)}) \cong
H^2(C^\cdot_{f_*})$. Thus, the obstructions lie  in the kernel of
$\mathcal{I} : H^2(C^\cdot_{f_*}) \lrg H^2(C^\cdot_{(\eta,\mu)})$.

As to $H^2(C^\cdot_{(\eta,\mu)})$, consider $K$ as in
Equation~(\ref{equa diagramma tre piani}), i.e.,
$$
K= \left\{f\in \Htp\left(\ker(\de),\dfrac{A_Z}{\de A_Z}\right)\mid
f(I_\Gamma\cap\ker(\de))\subset \dfrac{I_\Gamma}{I_\Gamma\cap\de
A_Z}\right\}
$$
and the exact sequence
$$
0 \lrg K \stackrel{\eta}{\lrg} \Htp\left(\ker
(\de),\dfrac{A_Z}{\de A_Z}\right)
\stackrel{\pi'}{\lrg}\coker(\eta)\lrg 0,
$$
with $\coker(\eta)= \Htp\left(I_\Gamma\cap \ker(\de)
,\dfrac{A_\Gamma}{ \de A_\Gamma}\right)$. Then, $
H^2(C^\cdot_{(\eta,\mu)}) \cong H^2(C^\cdot_{\pi' \circ \mu}) $.
Let $J$ be as in (\ref{equa diagramma tre piani}), i.e.,
$$
J= \left\{f\in \Htp\left(\ker(\de),\dfrac{A_Z}{\de A_Z}\right)\mid
f(\ker(\de)\cap q^* A_Y)=0 \right\};
$$
thus,
$$
\pi' \circ \mu: J \lrg \Htp\left(I_\Gamma\cap \ker(\de)
,\dfrac{A_\Gamma}{ \de A_\Gamma}\right),
$$
with
$$
\coker(\pi' \circ \mu)=\Htp(I_\Gamma \cap \ker (\de)\cap q^* A_Y,
\dfrac{A_\Gamma}{\de A_\Gamma}).
$$
Consider the map $\mathcal{I\, '}:H^2(C^\cdot_{\pi' \circ \mu})
\lrg H^1(\coker(\pi' \circ \mu))= H^1(\Htp(I_\Gamma \cap \ker
(\de)\cap q^* A_Y, \dfrac{A_\Gamma}{\de A_\Gamma}) )$. Since the
complex $\de A_\Gamma$ is acyclic, the projection
$$
 \Htp(I_\Gamma \cap \ker (\de)\cap q^* A_Y,\dfrac{A_\Gamma}{\de
A_\Gamma} ) \lrg  \Htp(I_\Gamma \cap \ker (\de)\cap q^* A_Y,
A_\Gamma)
$$
is a quasi-isomorphism.

Therefore, the obstruction space is contained in the kernel of
$\mathcal{J}: H^2(C^\cdot_{f_*}) \lrg  H^1(\Htp(I_\Gamma \cap
\ker(\de)\cap q^*A_Y,A_\Gamma))$, i.e.,
\begin{center}
$\xymatrix{ H^2(C^\cdot_{f_*}) \ar[rrd]^{\mathcal{J}}
\ar[r]^{\mathcal{I}}  & H^2(C^\cdot_{\pi' \circ
\mu})\ar[r]^{\mathcal{I\, '}\qquad \qquad\ \  } &
H^1\left(\Htp(I_\Gamma \cap \ker (\de)\cap q^* A_Y,
\dfrac{A_\Gamma}{\de A_\Gamma}) \right)\ar[d]^\cong  \\
& & \oplus_i \Hom(H^i(I_\Gamma \cap \ker (\de)\cap q^* A_Y),
H^i(A_\Gamma)).
 \\ }$
\end{center}

\end{proof}

\begin{cor}\label{coroll semiregolari-ostruzion}
Let $f:X \lrg Y$ be a  holomorphic map of  compact K\"{a}hler
manifolds. Let $p=dim Y - dim X$. Then, the obstruction space to
the infinitesimal deformations of $f$  with fixed $Y$  is
contained in the kernel of the map
$$
\sigma : H^2(C^\cdot_{f_*}) \lrg H^{p+1}(Y,\Omega_Y^{p-1}).
$$
\end{cor}

\begin{proof}
Let $n=dim X$, $p=dim Y - dim X$  and $\mathcal{H}$ be the space
of harmonic forms on $Y$ of type $(n+1,n-1)$.  Using the
contraction with the forms $\omega \in  \mathcal{H}$, we define
the semiregularity map $\sigma$ as in Section~\ref{sezio def
semiregular map s}. Since $f^* \omega=0$, Lemma~\ref{lemma
f^*(x-|w)=f_*x-|f^*w} implies that the diagram
\begin{center}
$\xymatrix{A_X^{0,*}(f^*\Theta_Y)  \ar[r]^{\ \ \contr \omega}
& A_X^{n,*+n-1} \\
A_X^{0,*}(\Theta_X) \ar[u]^{f_*} \ar[r]  & 0 \ar[u]_\alpha   \\
}$
\end{center}
is commutative and we get the semiregularity map
$$
\sigma :H^2(C^\cdot_{f_*}) \lrg H^{p+1}(Y,\Omega_Y^{p-1}).
$$
Since $ q^* \mathcal{H}$ is contained in $I_\Gamma \cap \ker
\debar \cap \ker \de \cap q^* A_Y$, we conclude the proof applying
Theorem~\ref{teo ostruzi in ker H^2}.

\end{proof}

\begin{oss}

As we already noticed in Remark \ref{oss de-debar lemma per
X,Y,XxY, Gamma}, the previous corollary holds if the compact
complex manifolds are bimeromorphic to K\"{a}hler manifolds.

\end{oss}


\begin{thebibliography}{99}

\bibitem{bib Artin} M. Artin,
{\sl Deformations of Singularities,} Tata Institute of
Foundamental Research, Bombay, (1976).




\bibitem{bib behr-fante} K. Behrend, B. Fantechi,
\emph{The Intrinsic Normal Cone}, Invent. Math., {\bf 128} (No.
1), (1997), 45-88.


\bibitem{bib bloch} S. Bloch,
\emph{Semi-regularity and de Rham cohomology}, Invent. Math., {\bf
17}, (1972), 51-66.


\bibitem{bib buchFlenner} R-O. Buchweitz, H. Flenner,
\emph{A semiregularity map for modules and applications to
deformations}, Compos. Math., {\bf 137} (No. 2), (2003), 135-210.

\bibitem{bib buchMilso} R-O. Buchweitz, J.J. Milson,
{\sl CR-geometry and deformations of isolated singularities}, Mem.
Amer. Math. Soc.  {\bf 125} (1997), no. 597.


\bibitem{bib clemens} H. Clemens,
\emph{Geometry of formal Kuranishi theory,}  Adv. Math., {\bf 198}
(No. 1), (2005),   311-365.


\bibitem{bib DeligneGMS} P. Deligne, P. Griffiths, J. Morgan,
 D. Sullivan,
 \emph{Real
homotopy theory of K¨ahler manifolds,} Invent. Math. {\bf 29}
(1975), 245-274.


\bibitem{bib Fantec-Manet}B. Fantechi, M. Manetti,
\emph{Obstruction calculus for functors of Artin rings, I,} J.
Algebra, {\bf 202}, (1998), 541-576.

\bibitem{bib Fantec-Manet 2} B. Fantechi, M. Manetti,
\emph{On the $T^1$-lifting theorem,} J.~Algebraic Geom., {\bf 8},
(1999), 31-39.


\bibitem{bib fiorenz-Manet} D. Fiorenza, M. Manetti,
\emph{$L_{\infty}$ algebras, Cartan homotopies and period maps,} Preprint
\texttt{arXiv:math.AG/0605297v1}.







\bibitem{bib HorikawaI} E. Horikawa,
\emph{On deformations of holomorphic maps I, II, } J. Math. Soc.
Japan, {\bf 25} (No.3), (1973), 372-396; {\bf 26} (No.4), (1974),
647-667.

%\bibitem{bib Horikawa III} E. Horikawa,
%\emph{On deformations of holomorphic maps III, }   Math. Annalen,
%{\bf 222}, (1976), 275-282.


\bibitem{bib tesidottorato} D. Iacono,
{\sl Differential Graded Lie Algebras and Deformations of
Holomorphic Maps,} Phd. Thesis, Roma, (2006),
\texttt{arXiv:math.AG/0701091}


\bibitem{bib articolo1tesi} D. Iacono,
{\sl $L_\infty$-algebras and deformations of Holomorphic Maps,}
Int. Math. Res. Not.~\textbf{8} (2008) 36 pp;
\texttt{arXiv:0705.4532v2}.



\bibitem{bib Kawamata}  Y.~Kawamata,
\emph{Unobstructed deformations, II}, J.~Algebraic Geom.,
\textbf{4}, (1995), 277-279.



\bibitem{bib kontsevich} M. Kontsevich,
\emph{Deformation quantization of Poisson manifolds, I.} Letters
in Mathematical Physics, {\bf 66}, (2003), 157-216;
\texttt{arXiv:q-alg/9709040}.


\bibitem{bib mane COSTRAINT} M. Manetti,
\emph{Cohomological constraint to deformations of compact K\"ahler
manifolds.}, Adv. Math.,  {\bf 186}, (2004), 125-142;
\texttt{arXiv:math.AG/0105175}.

\bibitem{bib manRENDICONTi} M. Manetti,
\emph{Lectures on deformations of complex manifolds,} Rend. Mat.
Appl. (7), {\bf24}, (2004), 1-183; \texttt{arXiv:math.AG/0507286}.

\bibitem{bib ManettiPREPRINT} M. Manetti,
\emph{Lie description of higher obstructions to deforming
submanifolds,} Ann. Sc. Norm. Super. Pisa Cl. Sci. \textbf{6}
(2007) 631-659;  \texttt{arXiv:math.AG/0507287}.




\bibitem{bib Ran ostruzio}Z. Ran,
\emph{Semiregularity, obstructions and deformations of Hodge
classes}, Ann. Scuola Norm. Pisa Cl. Sci., {\bf 28} (No. 4),
(1999), 809-821.

\bibitem{bib Ran ostruzioII}Z. Ran,
\emph{Universal variations of Hodge structure an
Calabi-Yau-Schottky relations}, Invent. Math., {\bf 138}, (1999),
425-449.



\end{thebibliography}
\end{document}